\documentclass[11pt]{article}
\usepackage{a4}
\usepackage{amsfonts}
\usepackage{latexsym}
\usepackage{url}

\begin{document}

\newcommand{\comment}[1]{}

\def\tuple{\ \hat{}\ }
\def\td#1{{\rm tr.deg}(#1)}

\def\o{\overline}

\def\R{\mathbb{R}}
\def\Q{\mathbb{Q}}
\def\N{\mathbb{N}}
\def\Z{\mathbb{Z}}
\def\C{\mathbb{C}}

\newtheorem{conjecture}{Conjecture}
\newtheorem{theorem}{Theorem}
\newtheorem{lemma}{Lemma}
\newtheorem{proposition}{Proposition}
\newtheorem{corollary}{Corollary}
\newtheorem{definition}{Definition}
\newtheorem{problem}{Problem}
\newtheorem{remark}{Remark}
\newtheorem{example}{Example}
\newtheorem{hypothesis}{Hypothesis}

\makeatletter
\def\@yproof[#1]{\@proof{ #1}}
\def\@proof#1{\begin{trivlist}\item[]{\em Proof#1.}}
\newenvironment{proof}{\@ifnextchar[{\@yproof}{\@proof{} 
}}{~$\Box$\end{trivlist}}
\makeatother

\title{The Theory of Liouville Functions}

\author{Pascal Koiran\\
LIP, Ecole Normale Sup\'erieure de Lyon\\
46, all\'ee d'Italie\\
69364 Lyon Cedex 07, France\\
\url{Pascal.Koiran@ens-lyon.fr}\\
\url{http://www.ens-lyon.fr/~koiran}}

\maketitle

\begin{abstract}
A Liouville function is an analytic function $H: \C \rightarrow \C$ 
with a Taylor series $\sum_{n=1}^{\infty} x^n/a_n$ such 
the $a_n$'s form a ``very fast growing'' sequence of integers.
In this paper we exhibit the complete first-order theory of the
complex field expanded with $H$.
\end{abstract}

\section{Introduction}

In~\cite{Wilkie00} Wilkie calls ``Liouville function'' 
a function $H:\C \rightarrow \C$ with a Taylor series of the form 
$$H(x)=\sum_{i=1}^{\infty} x^i/a_i$$
where the $a_i$ are non-zero integers satisfying the condition:
\begin{equation} \label{growth}
\mbox{for every } l \geq 1,\ |a_{i+1}|  > |a_i|^{i^l} 
\mbox{for all sufficiently large } i.
\end{equation}
A fragment of the first-order theory of the complex field expanded
with $H$ is described in~\cite{Wilkie00}.
In this paper we exhibit the complete first-order theory.
It turns out that this theory is the ``limit theory of generic
polynomials'' recently studied in~\cite{Koi01a} (this answers a
question of Zilber~\cite{Zil01a}).
We recall the axiomatization of the theory in section~\ref{axioms}, 
where we also present another equivalent axiomatization which is
closer in spirit to~\cite{Wilkie00}.
In section~\ref{structure} we give a short overview of the proof of
our main result and present two of the main tools: continuity of the
roots of polynomial systems and effective quantifier elimination.
The last two sections are devoted to the proof of the main result.

A  model of our theory can be constructed by a
Hrushovski-style amalgamation method~\cite{Koi01a,Zilgen}.
It is therefore natural to ask whether analytic models exist for other
theories constructed by this method.
The limit theory of generic
curves~\cite{CHKP} and the theory of generic functions with
derivatives~\cite{Zilgen} are two natural candidates. 
Additional examples and further discussion can be found in the
surveys~\cite{Poi00,Zil01a}.

\section{Axiomatization} \label{axioms}

In this section only we work within an arbitrary algebraically closed
field $K$ of characteristic 0. In the remainder of the paper we set
$K=\C$. We do not work in the language of ``curved fields'' (the language of
fields expanded with a binary predicate) as in~\cite{Koi01a} but in the
language of fields expanded with a unary function symbol $H$.
We call this language $\cal L$.
Given a tuple $\o{x}$ of elements of $K$, 
$H(\o{x})$ denotes the tuple obtained by applying $H$ componentwise.
This notation will be used freely throughout the paper, for $H$ as
well as for other unary functions.

\subsection{The Limit Theory of Generic Polynomials} \label{generic}

A generic polynomial of degree $d$ is of the form
$g_d(x)=\sum_{i=1}^d \alpha_i x^i$ where the coefficients $\alpha_i$
are algebraically independent over $\Q$.
Let $F$ be a sentence of~$\cal L$.
We have shown in~\cite{Koi01a} that $F$ is either true for all generic
polynomials of sufficiently high degree, or false for all generic
polynomials of sufficiently high degree.
The set $T$ of sentences which are ultimately true therefore forms a
complete theory.
We recall  that this theory is defined by the
following axioms.
\begin{enumerate}
\item The axioms of algebraically closed fields of characteristic 0.
\item $H(0)=0$.
\item The universal axioms. Let $\phi(x_1,\ldots,x_n,y_1,\ldots,y_n)$ be a
conjunction of polynomial equations with coefficients in~$\Q$. 
If the subset of $K^{2n}$ defined by $\phi$ is of dimension
$<n$, we add the axiom 
\begin{equation} \label{aaxiom}
\forall x_1,\ldots,x_n\ \bigwedge_{i} x_i \neq 0 \wedge 
\bigwedge_{i \neq j} x_i \neq x_j 
\rightarrow \neg \phi(\overline{x},H(\o{x})).
\end{equation}

\item The inductive axioms. Let $\phi(x_1,y_1,\ldots,x_n,y_n,\overline{z})$
 be a conjunction of polynomial equations with rational coefficients.
For any fixed value of the parameter $\overline{z}$, $\phi$ defines an
algebraic subset $V_{\overline{z}}$ of $K^{2n}$.
Let $\xi(\overline{z})$ 
be a formula of the language of fields which states that  
$V_{\overline{z}}$ 
 is irreducible, has dimension $n$ 
 and is not contained in a subspace of the form
 $x_i = x_j$ for some $i \neq j$, or of the form $x_i  = c$ for some
element $c$ in the model.

Let $\epsilon$ be a function which chooses one variable 
$u_i^{\epsilon} \in \{x_i,y_i\}$ for every $i \in \{1,\ldots,n\}$.
For each value of the parameter $\overline{z}$,
the formula $\exists u_1^{\epsilon},\ldots,u_n^{\epsilon} 
\phi(\overline{x},\overline{y},\overline{z})$ defines a constructible
set $C^{\epsilon}_{\overline{z}} \subseteq K^n$.
As pointed out in ~\cite{CHKP}, there is a formula 
$\psi_{\epsilon} (\overline{z})$ 
of the language of  fields which states that
$C^{\epsilon}_{\overline{z}}$ is dense in $K^n$.
Let $\psi(\overline{z})$ be the disjunction of the $2^n$ formulas 
$\psi_{\epsilon} (\overline{z})$.
Let $\theta$ be the conjunction of $\xi$ and $\psi$.
We add the following axiom:
\begin{equation} \label{aeaxiom}
\forall \overline{z}\ \exists x_1,\ldots,x_n\
\theta(\overline{z}) 
\rightarrow  \phi(\overline{x},H(\overline{x}),\overline{z}).
\end{equation}
\end{enumerate}
In fact these inductive axioms are slightly different from those
of~\cite{Koi01a}. Indeed, there was not requirement that $\dim V_{\o{z}}=n$ 
in that paper. 
An inspection of the completeness proof in~\cite{Koi01a} reveals that
the inductive axioms are used only in the case $\dim V_{\o{z}}=n$.
The corresponding theories are therefore identical.

We have shown in~\cite{Koi01a} that a universal axiom is satisfied
by a generic polynomial $g_d$ of degree $d$ as soon as $d \geq n$ 
and that an inductive axiom is satisfied as
soon as $d \geq n(nr+n+r+1)$.
More precisely, we have obtained the following result.
\begin{theorem} \label{inductive}
Let $\phi(x_1,\ldots,x_n,y_1,\ldots,y_n,\o{z})$ be a
conjunction of polynomial equations with coefficients in $\Q$.
Given a tuple $\o{z}$ of parameters 
we denote by $V_{\o{z}}$ the algebraic subset of $\C^{2n}$ defined by
$\phi$.

Fix a tuple $\o{z}$ which satisfies the associated condition $\theta(\o{z})$
from~(\ref{aeaxiom}). 
Let $r$ be the transcendence degree of $\Q(\o{z})$ over $\Q$.
For any $k \geq 0$, if $d \geq n(nr+n+r+1)$ 
there exists $\o{x} \in \C^n$ such that $(\o{x},g_d(\o{x}))$ is a
generic point of $V_{\o{z}}$ (i.e., a point of 
transcendence degree $n$ over $\Q(\o{z})$).
\end{theorem} 

\subsection{Axiomatization \`a la Wilkie} \label{wilkie}

Given $I \subseteq \{1,\ldots,n\}$, say $I=\{i_1,\ldots,i_k\}$ (in
increasing order), and a $n$-tuple $x=(x_1,\ldots,x_n)$ of elements or
variables, we denote by
$\o{x}_I$ the $k$-tuple $(x_{i_1},\ldots,x_{i_k})$.

Consider $n$ polynomials
$f_i(x_1,\ldots,x_n,y_1,\ldots,y_n,z_1,\ldots,z_r)$ with integer
coefficients.
For each value of the parameter $\o{z}$ we obtain a polynomial
map $F_{\o{z}}: K^n \times K^n \rightarrow K^n$.
Recall that a zero $(\o{a},\o{b})$ of $F_{\o{z}}$ is said to be 
regular if the Jacobian matrix of $F_{\o{z}}$ at  $(\o{a},\o{b})$ has
rank $n$, or in other words if
\begin{equation} \label{det}
\det \left( {\partial F_{\o{z}} \over \partial
(\o{x}_I,\o{y}_J)}(\o{a},\o{b}) \right) \neq 0
\end{equation}
for some $I,J \subseteq \{1,\ldots,n\}$ such that $|I|+|J|=n$.
Let $V(F_{\o{z}})$ be the set of zeros of $F_{\o{z}}$.
A regular zero lies on a unique irreducible component of $V(F_{\o{z}})$
and this component is of dimension $n$.
Following Wilkie, we say that $(\o{a},\o{b})$ is a balanced
zero\footnote{Wilkie uses the terminology ``balanced, non-singular
zero''. We just write ``balanced zero'' for short.} of
$F_{\o{z}}$ if $a_1,\ldots,a_n$ are all non-zero and
pairwise distinct, and if one can choose $I$ and $J$ so that $I \cup J =
\{1,\ldots,n\}$ (or equivalently, so that $I \cap J = \emptyset$).

In a geometric language, condition~(\ref{det}) 
means that the tangent space to 
$V(F_{\o{z}})$ at $(\o{a},\o{b})$ has dimension $n$. 
If this tangent space is not included in a subspace of the form 
$x_i = c$ for some $i \in \{,\ldots,n\}$ and some constant $c \in K$
and if additionally  $(\o{a},\o{b})$ is a balanced zero of $F_{\o{z}}$,
we say that $(\o{a},\o{b})$ is a {\em well balanced zero}.

In this axiomatization we keep the axioms 1,2 and 3 from
section~\ref{generic} but we replace the inductive axioms by:
\begin{itemize}
\item[4'.] Let $\theta'(\o{z})$ be a formula of the language of fields which
expresses that $F_{\o{z}}$ has a well balanced zero.
We add the following axiom.

\begin{equation} \label{balaxiom}
\forall z_1,\ldots,z_r\ \exists x_1,\ldots,x_n\
\theta'(\o{z}) \rightarrow F(\o{x},H(\o{x}),\o{z})=0.
\end{equation}
\end{itemize}

The following lemma is standard. 
\begin{lemma} \label{regular}
Let $(\o{a},\o{b})$ be a regular zero of $F_{\o{z}}$.
For any neighbourhood $V$ of $(\o{a},\o{b})$ in $\C^{2n}$ there exists a
neighbourhood $U$ of $\o{z}$ such that for any $\o{\zeta} \in U$,
$F_{\o{\zeta}}$ has a regular zero in $V$.
\end{lemma}
\begin{proof}
Let $I$ and $J$ be such that condition~(\ref{det}) is satisfied.
There exist $n$ affine functions $l_1,\ldots,l_n$ such that
$(\o{a},\o{b})$ is an isolated solution of the system
$f_1=0,\ldots,f_n=0,l_1=0,\ldots,l_n=0$. 
Since we now have as many equations as unknowns (namely, $2n$)
we can apply Proposition~\ref{continuity}: there exists a
neighbourhood $U$ of $\o{z}$ such that for any $\o{\zeta} \in U$,
the system $F_{\o{\zeta}}=0,l_1=0,\ldots,l_n=0$ has a zero
$(\o{\alpha},\o{\beta}) \in V$.
Since $\det \left( {\partial F_{\o{\zeta}} \over \partial
(\o{x}_I,\o{y}_J)}(\o{\alpha},\o{\beta}) \right)$ 
is a continuous function of $\o{\zeta}$, $\overline{\alpha}$ and
$\o{\beta}$ we can choose $U$ so small that 
$\det \left( {\partial F_{\o{\zeta}} \over \partial
(\o{x}_I,\o{y}_J)}(\o{\alpha},\o{\beta}) \right) \neq 0$.
\end{proof}
The point of working with this second axiomatization 
is that we have the
following proposition.
\begin{proposition} \label{open}
The set of parameters $\o{z}$ such that $F_{\o{z}}$ has a balanced
zero is an open subset of $\C^r$.
The same is true of the set of parameters $\o{z}$ such that
$F_{\o{z}}$ has a well balanced zero.
\end{proposition}
\begin{proof}
Note that in the proof of Lemma~\ref{regular}, 
the subsets $I,J \subseteq \{1,\ldots,n\}$ which witness the fact that 
$F_{\zeta}$ has a regular zero are the same for all $\zeta \in U$.
This implies immediately the first part of the proposition.
The second  part follows from a similar continuity argument.
\end{proof}

\subsection{Equivalence of these axiomatizations} \label{equivalence}

We first show that any model of the limit theory of generic curves is
also a model of the theory defined in section~\ref{wilkie}.

Let $(K,H)$ be a model of the limit theory of generic polynomials.
Let $F: K^n \times K^n \times K^r \rightarrow K^n$ be a polynomial
map with integer coefficients.
Fix a tuple $\o{z}$ such that $F_{\o{z}}$ has a well-balanced zero
$(\o{x},\o{y})$.
This well balanced zero lies on an irreducible component $V$ of
$V(F_{\o{z}})$ defined by a conjunction 
$\phi(\o{x},\o{y},\o{\zeta})$ of polynomial equations where the 
parameters $\o{\zeta}$ lie in the algebraic closure of $\Q(\o{z})$.
We claim that $\o{\zeta}$ satisfies the associated condition 
$\theta(\o{\zeta})$ from the inductive axioms.
Indeed, it follows from the implicit function theorem that 
$K \models \psi(\o{\zeta})$. 
Moreover, $V$ is not included in a subspace of the form $x_i = x_j$
for some $i \neq j$ since the components of $\o{x}$ are pairwise
distinct. Finally, $V$ is not included in a subspace of the form
$x_i=c$ for some constant $c$ due to the condition on the tangent
space at $(\o{x},\o{y})$.
We can therefore apply the inductive axioms: 
there exists $\o{x} \in \C^n$ such that $\phi(\o{x},H(\o{x}),\o{z})$.
In particular we have $F(\o{x},H(\o{x}),\o{z})=0$.

Next we show that any model of the theory defined in
section~\ref{wilkie} is a model of the limit theory of generic
polynomials. 
Consider therefore a model $(K,H)$ of the theory defined in
section~\ref{wilkie}. 
We assume without loss of generality that $K$ is of infinite
transcendence degree over $\Q$.
\begin{lemma} \label{trianglemma}
Let $S \subseteq K^m$ be a constructible set 
defined by a boolean combination of polynomial equations with 
coefficients in a subfield $k$ of $K$. 
Assume that the projection of $S$ on the first $d$
variables is dense in $K^d$. 

There exists a polynomial map $F = (f_1,\ldots,f_{m-d}):
K^m \rightarrow K^{m-d}$ 
which satisfies the following properties:
\begin{itemize}
\item[(i)] $f_i$ depends only on the
first $d+i$ variables and its coefficients are in $k$  
(i.e., $f_i \in k[X_1,\ldots,X_{d+i}]$).
\item[(ii)] The algebraic set $V(F) = \{x \in K^m;\ F(x) =0\}$ has
a dense projection on the first $d$ variables.

\item[(iii)] There exists a nonzero polynomial 
$P \in k[X_1,\ldots,X_d]$ such that for any $x \in V(F)$, 
$P(x_1,\ldots,x_d) \neq 0$ implies that $x \in S$ and that $x$ is a
regular zero of $F$ (more precisely, the matrix of partial derivatives
of $f_i$, $i=1,\ldots,m-d$ with respect to $X_j$, $j=d+1,\ldots,m$
has rank $m-d$ at $x$).
\end{itemize}
\end{lemma}
Note that property (iii) implies (ii) by the implicit function theorem. 
We will use this lemma only in the case where $S$ is an algebraic set.
\begin{proof}[of Lemma~\ref{trianglemma}]
We assume without loss of generality that $k$ is of finite
transcendence degree over $\Q$.
Since the projection of $S$ on the first $d$ variables is dense,
there exists a point $(\alpha_1,\ldots,\alpha_d)$ of transcendence $d$
over $k$ which is in the projection.
In fact, there exists $\alpha \in S$ such that $\alpha_1,\ldots,\alpha_d$ 
is a transcendence basis of $\alpha$ over $k$.
Let $f_i$ be the minimal polynomial of $\alpha_{d+i}$ over 
$k(\alpha_1,\ldots,\alpha_{d+i-1})$.
Condition~(i) is satisfied by definition, and condition (ii) is also
satisfied since $\alpha \in V(F)$.

Let $(x_{d+1},\ldots,x_m)$ be such that 
$(\alpha_1,\ldots,\alpha_d,x_{d+1},\ldots,x_m) \in V(F)$
and let $\phi(X_1,\ldots,X_m)$  be a a boolean combination 
of polynomial equations with coefficients in $k$
which is satisfied by $\alpha_1,\ldots,\alpha_m$.
This formula is also satisfied by 
$(\alpha_1,\ldots,\alpha_d,x_{d+1},\ldots,x_m)$
since the fields $k(\alpha_1,\ldots,\alpha_d,x_{d+1},\ldots,x_m)$
and $k(\alpha_1,\ldots,\alpha_m)$ are isomorphic.
In other words, $(\alpha_1,\ldots,\alpha_d)$ 
satisfies formula $\Phi'(x_1,\ldots,x_d)$ below:
$$\forall x_{d+1},\ldots,x_m F(x_1,\ldots,x_m) = 0
\Rightarrow \Phi(x_1,\ldots,x_m).$$ 
Since $\alpha_1,\ldots,\alpha_d$ are algebraically independent over $k$, 
$\Phi'$ defines a Zariski dense subset of $K^d$.
This is exactly the first part of condition (iii), 
if we take for $\Phi$ the formula defining $S$.
To obtain the second part of this condition, we apply the same
observation to a different $\Phi$. 
Namely, we apply it to the formula 
$\displaystyle \Phi(x_1,\ldots,x_m) \equiv \bigwedge_{i=1}^{m-d} 
{\partial f_i \over \partial X_{d+i}}(x_1,\ldots,x_{d+i}) \neq 0.$
This formula is satisfied by $\alpha$ since $f_i$ is 
the minimal polynomial of $\alpha_{d+i}$ over 
$k(\alpha_1,\ldots,\alpha_{d+i-1})$.
Our observation now implies that exists a nonzero polynomial 
$R \in k[X_1,\ldots,X_d]$ such that for any point $x \in V(F)$ with 
$R(x_1,\ldots,x_d) \neq 0$, the $m-d$ partial derivatives 
$\partial f_i \over \partial X_{d+i}$ do not vanish at $x$.
Since the jacobian matrix of $F$ contains a triangular matrix 
with these partial derivatives on the diagonal, it has maximum rank
$m-d$ at $x$ and this point is by definition a regular zero of $F$.
\end{proof}
There is of course nothing special about the first $d$ variables in
this lemma: we may project on any tuple of variables as long as the
projection is dense. 
This is just what we shall do now.

Let $\phi(x_1,y_1,\ldots,x_n,y_n,\overline{z})$ be a conjunction of
polynomial equations with rational coefficients.
Fix $\overline{z}$ such that the associated formula $\theta$
in~(\ref{aeaxiom}) is satisfied.
Since $K \models \theta(\overline{z})$ 
there exists $I,J \subseteq \{1,\ldots,n\}$ such that $|I|+|J|=n$,
$I \cap J = \emptyset$ and the projection 
of $V_{\o{z}}$ on the variables 
$\o{x}_I \tuple \o{y}_J$ is dense in $K^n$.
Let us apply Lemma~\ref{trianglemma} to $S=V_{\o{z}}$.
We obtain a polynomial map 
$F:K^n \times K^n \rightarrow K^n$ and a polynomial $P$ in $n$ variables 
with coefficients in $k=\Q(\o{z})$ 
such that for any point $\o{x}\tuple \o{y} \in V(F)$, 
$P(\o{x}_I,\o{y}_J) \neq 0$ 
implies that $\o{x}\tuple \o{y}  \in V_{\o{z}}$ 
and that $\o{x}\tuple \o{y}$ is a
regular zero of $F$.
Since $\dim V_{\o{z}} = n$ and $V_{\o{z}}$ is irreducible, 
this variety is an irreducible component of $V(F)$.
Now we consider the polynomial map $G:K^{2(n+1)} \rightarrow K^{n+1}$ 
which sends $(x_1,\ldots,x_n,x_{n+1},y_1,\ldots,y_n,y_{n+1})$ to
$F(x_1,\ldots,x_n,y_1,\ldots,y_n)\tuple  
f_{n+1}(x_1,\ldots,x_n,x_{n+1},y_1,\ldots,y_n,y_{n+1})$.
Here $x_{n+1}$ and $y_{n+1}$ are two additional variables, and 
$f_{n+1}
 = P(\o{x}_I,\o{y}_J)y_{n+1} -1$.
We claim that $G$ has a well balanced  zero.
To obtain such a zero, pick a generic point $\o{a}\tuple \o{b}$ of
$V_{\o{z}}$. 
Since $K \models \theta(\o{z})$, the components $a_1,\ldots,a_n$ are
all nonzero and distinct from each other.
Moreover we can set 
$b_{n+1} = 1/P(\o{a}_I,\o{b}_J)$ since $\o{a}_I\tuple \o{b}_J$ has
transcendence degree $n$ over $k$.
Pick an arbitrary $a_{n+1}$ different from 0
and from $a_1,\ldots,a_n$.
The matrix of partial derivatives of $f_1,\ldots,f_{n+1}$ at
$(a_1,\ldots,a_{n+1},b_1,\ldots,b_{n+1})$ 
with respect to the variables  $x_i$ ($i \in J$) 
and $y_i$ ($i \in I \cup \{n+1\}$)
has the block form
$$B = \left(\begin{array}{cc}
A & 0 \\
0 & P(\o{a}_I,\o{b}_J)
\end{array}
\right)$$
where $A$ is the matrix of partial derivatives of $f_1,\ldots,f_n$
with respect to the variables $x_i$ ($i \in J$) 
and $y_i$ ($i \in I$).
We know from Lemma~\ref{trianglemma} that $A$ has rank~$n$.
Hence $B$ has rank $n+1$ and $(a_1,\ldots,a_{n+1},b_1,\ldots,b_{n+1})$ 
is a balanced zero of $G$.
In order to show that this zero is well balanced, we have to check the
condition on the tangent space. 
This condition is indeed satisfied due to Lemma~\ref{vertical} below
and to the fact $x_{n+1}$ does not
appear in $f_1,\ldots,f_{n+1}$.
\begin{lemma} \label{vertical}
The tangent space to $V(F)$ at $\o{a}\tuple \o{b}$ is not included in a
subspace of the form $x_i=c$ for some constant $c \in K$.
\end{lemma}
\begin{proof}
Assume the opposite. We have seen that $V_{\o{z}}$ is an irreducible
component of $V(F)$.
Since $\o{a}\tuple \o{b}$ is a generic point of $V_{\o{z}}$ this
variety would be included in the subspace $x_i=c$. 
This is in contradiction with the hypothesis $K \models \theta(\o{z})$.
\end{proof}
We can therefore apply the modified inductive axioms~(\ref{balaxiom}) to $G$: 
there exists $x_1,\ldots,x_{n+1}$ such that 
$F(x_1,\ldots,x_n,H(x_1),\ldots,H(x_n),\o{z})=0$
and $P(\o{x}_I,\o{y}_J)y_{n+1} -1=0$.
Since $P(\o{x}_I,\o{y}_J) \neq 0$, 
$(x_1,\ldots,x_n,H(x_1),\ldots,H(x_n)) \in V_{\o{z}}$ by
property~(iii): 
we have proved that the inductive axioms~(\ref{aeaxiom}) are satisfied.

\section{Structure of the Proof} \label{structure}

Wilkie has shown in~\cite{Wilkie00} 
that Liouville functions satisfy the universal axioms. 
In order to show that the theory of Liouville functions is the limit
theory of generic polynomials, it remains to show that the inductive
axioms are also satisfied.
Using the axiomatization of section~\ref{wilkie},
we need the following result.
\begin{theorem}[Main Theorem] \label{main}
Let  $G:\C^n \times \C^n \rightarrow \C^n$
be a polynomial map.
If $G$  has a well balanced
zero, there exists $\o{a} \in \C^n$ such that 
$G(\o{a},H(\o{a}))=0$.
\end{theorem}
The proof will be given at the end of section~\ref{path}.
A special case of this theorem (which we will not use here) was
obtained in~\cite{Wilkie00}: 
let $G:\C^n \times \C^n \rightarrow \C^n$ be a polynomial map with
integer coefficients.
If $G$ has a balanced zero then there  exists $\o{a} \in \C^n$ such that 
$G(\o{a},H(\o{a}))=0$.
Wilkie's proof of this result relies in particular on Newton's
method.\footnote{In fact he even shows that the map $\o{x} \mapsto
G(\o{x},H(\o{x}))$  has a non-singular zero with
pairwise distinct, nonzero coordinates.}
By contrast the proof of Theorem~\ref{main} is based on a method which
is reminiscent of homotopy methods for solving systems of polynomial
equations. Let $H_d(x)$ be the partial sum $\sum_{i=1}^d x^i / a_i$.
In section~\ref{startingpoint} we show that the system 
$G(\o{x},H_d(\o{x}))=0$ has an isolated solution 
$\o{x}_d \in \C^n$ for all sufficiently large $d$.
Then we track $\o{x}_d$ as $d$ goes to infinity.
It turns out that these roots remain in a compact subset of $\C^n$. 
Theorem~\ref{main} then follows immediately from a standard uniform
convergence argument. 

As in~\cite{Wilkie00} effective quantifier elimination also plays an
important role. Here there is an additional complication due to the
presence of arbitrary complex parameters in $G$.
Our solution to this problem is to work not for a single value of
the parameters (i.e. for a single $G$), 
but simultaneously for all parameters in a compact set.
This is made possible in particular by Proposition~\ref{open}.
We will use the following version of quantifier elimination (note that
we need to work over the real numbers).
\begin{proposition} \label{elimination}
Let $\Phi(x_1,\ldots,x_n)$ be a formula of the language of ordered
rings. The $k$ polynomials occuring in $\Phi$ have integer coefficients.
Let $h \geq 2$ be an upper bound on their absolute values.
Let $d \geq 2$ be an upper bound on the degrees of these polynomials,
and let $m$ be the number of occurrences of quantifiers in $F$.

In the theory of real-closed fields, $\phi$ is equivalent to a
quantifier-free formula $\Psi(x_1,\ldots,x_n)$ in which all polynomials
are of degree at most $d^c$, and have integer coefficients of absolute
value bounded by $h^{d^c}$. The constant $c$ depends only on $n$, $m$
and $k$.
\end{proposition}
Almost any reasonable quantifier elimination method will yield the
above result. Much more precise bounds are known, see for
instance~\cite{BaPoRoy96,Ren92}. 
We will not need them here since the parameters $n$, $m$ and $k$ can be
treated as constants for our purposes.

In the remainder of this section we present another important tool:
continuity of the roots of polynomial systems.
For $\o{z} = (z_1,\ldots,z_n) \in \C^n$ 
we set $||\o{z}||=\sum_{i=1}^n |z_i|^2$ 
(this is just the Euclidean norm on $\R^{2n}$).

\begin{proposition}[Continuity of roots] \label{continuity}
The following property holds for any polynomial map $F: \C^n
\times \C^r \rightarrow \C^n$ and any $\o{z} \in \C^r$.
 
Let $\o{x}$ be an isolated root of the map
$F_{\o{z}}: \o{x} \mapsto F(\o{x},\o{z})$.
For any sufficiently small neighbourhood $U$ of $\o{x}$ there exists a
neighbourhood $V$ of $\o{z}$ such that for all $\o{\zeta} \in V$,
the number of roots of $F_{\o{\zeta}}$ in $U$ is positive and finite.
\end{proposition}
This follows for instance from 
the ``extended geometric version'' of B\'ezout's
theorem~\cite{BCSS}.

\begin{corollary} \label{norm}
Let $F:\C^n \times \C^r \rightarrow \C^n$ be a polynomial map. 
The map 
$$\begin{array}{rcl}
N_F & : & \C \rightarrow \R \cup \{+\infty\}\\
& & z \mapsto \min \{||\o{x}||;\ \o{x} \mbox{ \rm is an isolated root of }
F(\o{x},\o{z})=0\}
\end{array}$$
is upper semi-continuous (we set $N_F(\o{z})=+\infty$ if the system
$F(\o{x},\o{z})=0$ has no isolated roots).
\end{corollary}
\begin{proof}
Fix $\o{z} \in \C^r$ such that $N_F(\o{z}) < + \infty$. 
We have to show that for every $\epsilon>0$ there is a neighbourhood
$U$ of $\o{z}$ such that $N_F(\o{z}) \geq N_F(\o{\zeta}) - \epsilon$
for every $\o{\zeta} \in U$. let $\o{x}$ be an isolated root of the
system $F(\o{x},\o{z})=0$ such that $N_F(\o{z})=||\o{x}||$.
By Proposition~\ref{continuity}.(i) there is a neighbourhood $U$ of
$\o{z}$ such that for every $\o{\zeta} \in U$ the system
$F(\o{x},\o{z})=0$ has an isolated root in the ball $B(x,\epsilon)$.
Hence $N_F(\o{\zeta}) \leq ||\o{x}|| + \epsilon = N_F(\o{z}) + \epsilon$.
\end{proof}
One can show that if $\o{z}$ is such that 
$F_{\o{z}}$ has finitely many roots then $N_F$ is
continuous in $\o{z}$.
The example of the polynomial $F(x,z)=zx^2-2x+1$ shows that 
 no such 
continuity property holds for the map 
$$z \mapsto \max \{||\o{x}||;\ \o{x} \mbox{ \rm is an isolated root of
} F(\o{x},\o{z})=0\}.$$
Indeed, $F$ has a single root for $z=0$; 
for $z {\not \in} \{0,1\}$ it has a second root which goes to infinity
as $z$ goes to 0.

\section{A Starting Point for the Homotopy} \label{startingpoint}

We denote $H_{d,\epsilon}(x) = \sum_{i=1}^d x^i/a_i + \epsilon
 x^{d+1}$ and $H_d = H_{d,0}$.
 Let $g_d$ be a generic polynomial of degree $d$.
The {\em modified partial sum} $\mu_{k,d}$ is the polynomial function
$x \mapsto H_k(x) + x^{k}g_d(x).$

\subsection{Finiteness for modified partial sums} \label{finitesec}

We temporarily revert to the language of curved fields to cite a
simple combinatorial result (Lemma~7 from~\cite{KP01a}).
\begin{lemma} \label{sufficient}
Let $k$ be a $p$-sufficient substructure of a curved field $(K,C)$
and $\o{z}=(z_1,\ldots,z_r)$ a tuple of $r$ elements of $K\backslash k$.
Set $q = \lfloor \frac{p-r}{r+1} \rfloor$.
Let $j$ be the smallest integer such that there exists an extension
$\ell$ of $k(\o{z})$ satisfying $\td{\ell/k} \leq p-q(j+1)$ 
and $\delta(\ell:k) \leq j$ (note that $j$ always exists and is upper
bounded by~$r$).
Then $\ell$ is $q$-sufficient. 
\end{lemma}
Some explanations are in order. In this lemma $K$ is an arbitrary
field and the ``curve'' $C$ is an arbitrary subset of $K^2$.
The symbol $\delta(\ell:k)$ is defined by the formula 
$$\delta(\ell:k) = \td{\ell:k}-{\rm Card}(C \cap \ell^2 - C \cap k^2).$$
A subfield $k$ of $K$ is said to be $p$-sufficient if
$\delta(\ell:k)\geq 0$ for any subfield $\ell$ of $K$ which
contains $k$, and is of transcendence degree at most $p$ over $k$.
For instance, we have seen is section~\ref{generic} that the universal
axiom~(\ref{aaxiom}) is satisfied by a generic polynomial of degree
$d$ as soon as $d \geq n$.
This implies that $\Q$ is $d$-sufficient if we interpret $C$ by the
graph of a generic polynomial $g:\C \rightarrow \C$ of degree $d$.
\begin{theorem} \label{finite}
Let $\phi(x_1,\ldots,x_n,y_1,\ldots,y_n,z_1,\ldots,z_r)$ be a
boolean combination of polynomial equations with coefficients in $\Q$.
Fix a tuple $\o{z}$ such that the constructible subset
$D_{\o{z}}$ of $\C^{2n}$ defined by $\phi(.,.,\o{z})$ 
has dimension at most $n$.
Let $C_d$ be the graph of a generic polynomial of degree $d$.
If $d \geq n(r+1)+r$ the system 
\begin{equation} \label{fsys}
\bigwedge_{i} x_i \neq 0  \bigwedge_{i \neq j} x_i \neq x_j 
\wedge \bigwedge_{i} C_d(x_i,y_i) \wedge
\phi(\o{x},\o{y},\o{z})
\end{equation}
has at most finitely many solutions in $\C^{2n}$.
\end{theorem}

\begin{proof}
Since $\Q$ is $d$-sufficient, 
by choice of $d$ and Lemma~\ref{sufficient}
there exists a $n$-sufficient extension
$\ell$ of $\Q(\o{z})$ of transcendence degree at most $d-n$ over~$\Q$.
There are at most $d-n$ nonzero points on $C_d$ with both coordinates in
$\ell^2$. Moreover, outside $\ell^{2n}$ any solution of the system 
$$\bigwedge_{i \neq j} x_i \neq x_j 
\wedge \bigwedge_{i} C_d(x_i,y_i)$$
must be
of transcendence degree at least $n$ over $\ell$ (by choice of $\ell$).
We conclude that up to a finite set, all solutions of~(\ref{fsys}) are
of transcendence degree at least $n$ over $\Q(\o{z})$.

This implies that the subset $S$ of $D_{\o{z}}$ 
defined by~(\ref{fsys}) is finite.
Indeed, if $S$ is infinite
this (constructible) subset of $D_{\o{z}}$ 
must contain infinitely many non-generic
points of $D_{\o{z}}$ (i.e., points of transcendence degree
over $\Q(\o{z})$ smaller than $ \dim D_{\o{z}}$).
\end{proof}
We have the same property for modified partial sums.
\begin{corollary} \label{finitemu}
Let $\phi(x_1,\ldots,x_n,y_1,\ldots,y_n,z_1,\ldots,z_r)$ be a
boolean combination
 of polynomial equations with coefficients in $\Q$.
Fix a tuple $\o{z}$ such that the algebraic subset
$D_{\o{z}}$ of $\C^{2n}$ defined by $\phi(.,.,\o{z})$ 
has dimension at most~$n$.
If $d \geq n(r+1)+r$ the system 
$$\bigwedge_{i} x_i \neq 0 \wedge  \bigwedge_{i \neq j} x_i \neq x_j 
 \wedge
\phi(\o{x},\mu_{k,d}(\o{x}),\o{z})$$
has at most finitely many solutions in $\C^{2n}$.
\end{corollary}
\begin{proof}
Let $\o{x}$ be a solution of the system.
Note that $(\o{x},g_d(\o{x}))$ lies in the constructible subset of
$\C^{2n}$ (call it $C_{\o{z}}$):
$$\bigwedge_{i} x_i \neq 0 \wedge
\phi(x_1,\ldots,x_n,H_k(x_1)+x_1^k y_1,\ldots,H_k(x_n)+x_n^k y_n).$$
Let $P:\C^{2n} \rightarrow \C^{2n}$ be the polynomial map
$$(x_1,\ldots,x_n,y_1,\ldots,y_n) \mapsto 
(x_1,\ldots,x_n,H_k(x_1)+x_1^k y_1,H_k(x_n)+x_n^k y_n).$$
Since $P(C_{\o{z}}) \subseteq D_{\o{z}}$ and every point in
$D_{\o{z}}$ has finitely many preimages in $C_{\o{z}}$ we have 
$\dim C_{\o{z}} \leq \dim D_{\o{z}} \leq n$. 
By Theorem~\ref{finite} (applied to
$C_{\o{z}}$) we conclude that our system has finitely many solutions.
\end{proof}

\subsection{Existence for Modified Partial Sums} \label{existsec}

The only property of the coefficients of Liouville functions that will
be used in the next proposition is that they are rational numbers.
\begin{proposition} \label{existence}
Let $\phi(x_1,\ldots,x_n,y_1,\ldots,y_n,\o{z})$ be a
conjunction of polynomial equations with coefficients in $\Q$.
Given a tuple $\o{z}$ of parameters 
we denote by $V_{\o{z}}$ the algebraic subset of $\C^{2n}$ defined by
$\phi$.

Fix a tuple $\o{z}$ which satisfies the associated condition $\theta(\o{z})$
from~(\ref{aeaxiom}). 
Let $r$ be the transcendence degree of $\Q(\o{z})$ over $\Q$.
For any $k \geq 0$, if $d \geq n(nr+n+r+1)$ 
there exists $\o{x} \in \C^n$ such that $(\o{x},\mu_{k,d}(\o{x}))$ is a
generic point over $\Q(\o{z})$ of $V_{\o{z}}$ (note that this
genericity condition implies in particular that the components of
$\o{x}$ are nonzero and pairwise distinct).
\end{proposition}
\begin{proof}
We argue as in Corollary~\ref{finitemu}.
One would like to find a point $(\o{x},g_d(\o{x}))$ on the
algebraic subset $W_{\o{z}}$ of  $\C^{2n}$ defined by the formula 
$$\phi(x_1,\ldots,x_n,H_k(x_1)+x_1^k y_1,\ldots,H_k(x_n)+x_n^k y_n).$$
Pick a generic point $(\o{\alpha},\o{\beta})$ of $V_{\o{z}}$. 
Then $(\o{\alpha},\o{\gamma}) \in W_{\o{z}}$ where $\gamma_i=(\beta_i-H_k(\alpha_i))/\alpha_i^k$.
Note that $\alpha_1,\ldots,\alpha_n$ are pairwise distinct and do not belong to
the algebraic closure $K$ of $\Q(\o{z})$. 
Moreover there exists $I,J \subseteq \{1,\ldots,n\}$ such that
$|I|+|J|=n$, $I \cap J = \emptyset$ and $\o{\alpha}_I\tuple \o{\beta}_J$ is of
transcendence degree $n$ over $K$.
The tuple $\o{\alpha}_I\tuple \o{\gamma}_J$ is also 
of transcendence degree $n$ over $K$ since 
$K(\o{\alpha}_I,\o{\beta}_J)=K(\o{\alpha}_I,\o{\gamma}_J)$.
We can therefore apply the inductive axioms to the irreducible
component of $W_{\o{z}}$ which contains $(\o{\alpha},\o{\gamma})$.
More precisely, by Theorem~\ref{inductive} there exists a point 
$(\o{x},g_d(\o{x})) \in W_{\o{z}}$ of transcendence degree $n$ over $K$.
We conclude that $(\o{x},\mu_{k,d}(\o{x}))$ is a point of  $V_{\o{z}}$ 
of transcendence degree $n$ over $K$.
\end{proof}

\subsection{Isolated Solutions}

The results of sections~\ref{finitesec} and~\ref{existsec} can be
summarized as follows.
\begin{theorem} \label{isomu}
Let $F:\C^n \times \C^n \times \C^r \rightarrow \C^n$ be a polynomial
map with integer coefficients.
For any $k \geq 0$, any $d \geq  n(nr+n+r+1)$ 
and any $\o{z}$ such that $F_{\o{z}}$ has a well
balanced zero,
the system $F(\o{x},\mu_{k,d}(\o{x}),\o{z})=0$ has an isolated solution.
\end{theorem}
\begin{proof}
Fix $k \geq 0$, $d \geq  n(nr+n+r+1)$ 
and $\o{z}$ such that $F_{\o{z}}$ has a well
balanced zero. 
This well balanced zero lies on an irreducible component $V$ of
$V(F_{\o{z}})$ of dimension $n$  
and $V$ satisfies the inductive
axioms. More precisely, we have seen at the beginning of
section~\ref{equivalence} that $V$ is defined by a conjunction 
$\phi(\o{x},\o{y},\o{\zeta})$ of polynomial equations where the 
parameters $\o{\zeta}$ lie in the algebraic closure of $\Q(\o{z})$ and
satisfy the associated condition $\theta(\o{\zeta})$. 
By Proposition~\ref{existence} there exists $\o{x} \in \C^n$ with nonzero,
pairwise distinct coordinates such that 
$(\o{x},\mu_{k,d}(\o{x}))$ is  in~$V$. 
Since $\dim V =n$,  by
Corollary~\ref{finite} there are only finitely many such $\o{x}$.
To make sure that $\o{x}$ is an isolated solution of the system
$F(\o{x},\mu_{k,d}(\o{x}),\o{z})=0$ it is therefore sufficient to
satisfy the following requirement: $(\o{x},\mu_{k,d}(\o{x}))$ should
lie on no other irreducible component of $V(F_{\o{z}})$ but $V$
(indeed, other components might be of dimension $>n$).
This is possible thanks to the genericity condition in the
conclusion of Proposition~\ref{existence}.
\end{proof}

Any inductive axiom is eventually satisfied by $H_d$ if $d$ is
sufficiently large. More precisely:
\begin{theorem} \label{sums}
Let $F:\C^n \times \C^n \times \C^r \rightarrow \C^n$ be a polynomial
map.
If $d$ is sufficiently large and $|\epsilon| \leq 1/|a_{d+1}|$ 
the following property holds:
for any $\o{z}$ such that $F_{\o{z}}$ has a well balanced zero, 
the system $F(\o{x},H_{d,\epsilon}(\o{x}),\o{z})=0$ has an isolated
solution. 
\end{theorem}
Note that we do not rule out the the possibility that 
non-isolated solutions might also exist.
\begin{proof}[of Theorem~\ref{sums}]
Set $d_0=n(nr+n+r+1)$. For any $d \geq 1$ and $\o{\alpha} \in \R^{d_0}$, 
let $$\nu_{k,\o{\alpha}}(x) = 
H_k(x)+x^k\sum_{j=1}^{d_0} \alpha_j x^j.$$
One can easily write down a formula $\phi(\o{\alpha})$
of the language of ordered
fields which expresses the fact for any $\o{z}$ such that $F_{\o{z}}$
has a well balanced zero, the system  
$F(\o{x},\nu_{k,\o{\alpha}}(\o{x}),\o{z})=0$ has an isolated
solution in~$\C^n$ (of course this involves the separation of the real
and imaginary parts of variables such as $x_1,\ldots,x_n$, which range
over the complex numbers).
It follows from Theorem~\ref{isomu} that $\R \models \phi(\o{\alpha})$
if $\alpha_1,\ldots,\alpha_{d_0}$ are algebraically independent.
By Proposition~\ref{elimination}, $\phi(\o{\alpha})$ is equivalent to a
quantifier-free formula $\psi(\o{\alpha})$ involving polynomials of
degree $k^{O(1)}$ with integer coefficients of absolute value
$|a_k|^{k^{O(1)}}$ (the implied constants may depend on $F$ but not on $k$).
Any $\o{\alpha}$ which is not a root of any of these polynomials will
satisfy $\phi$. If $k$ is sufficiently large 
and $|\epsilon| \leq 1/|a_{k+d_0}|$, it  is indeed the case 
that $(1/a_{k+1},\ldots,1/a_{k+d_0-1},\epsilon)$ is not a
root of any of these polynomials.
This follows from the growth rate
condition~(\ref{growth}). 
We conclude that for any $\o{z}$
such that $F_{\o{z}}$ has a well balanced zero, the system
$F(\o{x},H_{k+d_0-1,\epsilon}(\o{x}),\o{z})=0$ has an isolated solution.
\end{proof}
Instead of $\phi(\o{\alpha})$, one could (less easily) write down a
formula $\phi'(\o{\alpha})$ of the language of fields 
such that $\C \models \phi'(\o{\alpha})$ iff the system
$F(\o{x},\nu_{k,\o{\alpha}}(\o{x}),\o{z})=0$ has an isolated
solution for any $\o{z}$ such that $F_{\o{z}}$
has a well balanced zero.
Allowing the order relation 
 just makes it easier to express the fact that there
is an isolated solution. 
By contrast order plays an essential role in the proof of the next result.
\begin{corollary} \label{startbound}
Let $F:\C^n \times \C^n \times \C^r \rightarrow \C^n$ be a polynomial
map and let $K$ be a rational ball of $\C^r$ such that $F_{\o{z}}$ has a
well balanced zero for all $\o{z} \in K$.
There is a constant $c>0$ such that the following property holds if
$d$ is sufficiently large:
 for any $\o{z} \in K$, the system $F(\o{x},H_d(\o{x}),\o{z})=0$
has an isolated solution which satisfies $||\o{x}|| \leq |a_d|^{d^c}-1/d$.
\end{corollary}
By {\em rational ball} of $\C^r$ we mean a closed ball 
$B(\o{z},R) = \{\zeta \in \C^r;\ ||\zeta - z|| \leq R\}$ 
such that the radius $R$ is a rational number and the real and
imaginary parts of $\zeta_1,\ldots,\zeta_r$ are also rational.
\begin{proof}[of Corollary~\ref{startbound}]
By Theorem~\ref{sums}, 
if $d$ is sufficiently large (say, $d \geq d_0$) the system 
$F(\o{x},H_{d}(\o{x}),\o{z})=0$ has isolated solutions 
for all $\o{z} \in K$.
Pick any $d \geq d_0$ and consider the map $N_F:K \rightarrow \C$
which sends $x$ 
to $N_F(x)=\inf \{||\o{x}||;\ F(\o{x},H_{d}(\o{x}),\o{z})=0\}$.
By Corollary~\ref{norm} this function is upper semi-continuous on $K$.
Since $K$ is compact, $N_F$ reaches its (finite) supremum $R(d)$ on $K$.
Translating the definition of $R(d)$ in first-order logic immediately
yields a formula $\phi(u)$ in the language of ordered
rings such that $\R \models \forall u\ (\phi(u) \leftrightarrow u=R(d))$.
Note that $\phi$ has only rational parameters since $K$ is a rational ball.
By elimination of quantifiers from $\phi$  we conclude that $R(d)$ is
a root of a polynomial of degree $d^{O(1)}$ with integer coefficients
bounded in absolute value by $|a_d|^{d^{O(1)}}$ (the implied constants
may depend on $K$ and $F$ but not on $d$). 
We conclude that $R(d) \leq |a_d|^{d^{\alpha}}$ for some constant
$\alpha$, so that $R(d) \leq |a_d|^{d^{\alpha+1}}-1/d$ if $d$ is
sufficiently large.
\end{proof}
This property is not only valid for $K$ a rational ball: one could
generalize to arbitrary compact sets (see Corollary~\ref{compact}).
One may wonder why we insist on a bound of the form $|a_d|^{d^c}-1/d$ 
in Corollary~\ref{startbound} instead of e.g. $|a_d|^{d^c}$.
The reason will become apparent in the next section.

\section{The Path Following Method} \label{path}

In the preceding section we have proved the existence of isolated
roots for systems of the form $F(\o{x},H_d(\o{x}),\o{z})=0$.
In this section we show that some of the roots stay inside a fixed
compact ball as $d$ goes to infinity.
The main theorem then follows easily.

\begin{lemma} \label{uniform}
Let $G: \C^n \times \C^r \times \C \rightarrow \C^n$ be a polynomial
map. We denote by $G_{\o{z},\epsilon}$ 
the map $\o{x} \mapsto G(\o{x},\o{z},\epsilon)$.
Let $K$ be a compact subset of $\C^r$ and $C$ a compact subset of
$\C^n$ such that for any $\o{z} \in K$, $G_{\o{z},0}$ has
isolated roots in the interior of $C$.
There exists $\delta > 0$ such that $G_{\o{z},\epsilon}$ has isolated roots in
$C$ if $|\epsilon| \leq \delta$ 
and $\o{z} \in K$. 
\end{lemma}
\begin{proof}
By Proposition~\ref{continuity}, for any  $\o{\zeta} \in K$ there exists
$\delta(\o{\zeta})>0$ such that $G_{\o{z},\epsilon}$ 
has isolated roots in $C$ if $|\epsilon| \leq \delta(\o{\zeta})$ and $\o{z}$
belongs to the open ball $B(\o{\zeta},\delta(\o{\zeta}))$.
Since $K$ is covered by the open balls
$B(\o{\zeta},\delta(\o{\zeta}))$, by compactness there exists a
finite cover of the form 
$B(\o{\zeta_1},\delta(\o{\zeta_1})),\ldots,
B(\o{\zeta_k},\delta(\o{\zeta_k}))$.
Set $\delta=\min(\delta(\o{\zeta_1}),\ldots,\delta(\o{\zeta_k}))$.
Now fix any $\o{z} \in K$ and $\epsilon$ such that 
$|\epsilon| \leq \delta$.
Since $\o{z} \in B(\o{\zeta_i},\delta)$ for some $i$, 
$G_{\o{z},\epsilon}$ has isolated roots in $C$ by choice of $\delta$.
\end{proof}

\begin{proposition} \label{induction}
Fix a polynomial map 
$F:\C^n \times \C^n \times \C^r \rightarrow \C^n$
with integer coefficients 
and a rational ball $K \subseteq \C^r$ such that $F_{\o{z}}$ 
has a well-balanced zero for all $\o{z} \in K$.
For any $c>0$ 
the following property holds if $d$ is sufficiently large.

Suppose that $R \leq |a_d|^{d^c}$ is an integer such that for all 
$\o{z} \in K$ the system $F(\o{x},H_d(\o{x}),\o{z})=0$ has an isolated
solution in the closed ball $\o{B}(0,R-1/d)$. 
If  $|\epsilon| \leq 1/|a_{d+1}|$, for any $\o{z} \in K$
the system $F(\o{x},H_{d,\epsilon}(\o{x}),\o{z})=0$
has an isolated solution in the closed ball $\o{B}(0,R-1/(d+1))$.
\end{proposition}
\begin{proof} 
By Theorem~\ref{sums}
if $d$ is sufficiently large (say, $d \geq d_0$) for any $\o{z} \in K$
the system  $F(\o{x},H_{d,\epsilon}(\o{x}),\o{z})=0$ has isolated solutions.
Pick any $d \geq d_0$ and 
let $R \leq |a_d|^{d^c}$ be an integer such that for all 
$\o{z} \in K$ the system $F(\o{x},H_d(\o{x}),\o{z})=0$ has an isolated
solution in the closed ball $\o{B}(0,R-1/d)$. 

Let $C = \o{B}(0,R-1/(d+1))$.
We can apply Lemma~\ref{uniform} to
$G(\o{x},\o{z},\epsilon)=F(\o{x},H_{d,\epsilon}(\o{x}),\o{z})$:
there exists $\delta>0$ such that for any
$\o{z} \in K$, the system $F(\o{x},H_{d,\epsilon}(\o{x}),\o{z})=0$ has
an isolated root $\o{x}$  in $C$ if $|\epsilon| \leq \delta$.
The same quantifier elimination argument as in
Corollary~\ref{startbound} shows that one may take
$1/\delta=|a_d|^{d^{O(1)}}$ (the implied constant depends only on $K$,
$F$ and $c$). One may therefore take $\delta=1/|a_{d+1}|$ if $d$ is
sufficiently large.
\end{proof}

\begin{theorem} \label{boundedroots}
Fix a polynomial map 
$F:\C^n \times \C^n \times \C^r \rightarrow \C^n$
with integer coefficients 
and a rational ball $K \subseteq \C^r$ such that $F_{\o{z}}$ 
has a well-balanced zero for all $\o{z} \in K$.
There exists $R>0$ such that the following property 
holds for all sufficiently large $d$: for any $\o{z} \in K$ the system 
$F(\o{x},H_d(\o{x}),\o{z})=0$ has an isolated
 root in the closed ball $\o{B}(0,R)$.
\end{theorem}
\begin{proof}
By Corollary~\ref{startbound} there is a constant $c>0$ 
such that if $d$ is sufficiently large (say, $d \geq d_0$)
then for any $\o{z} \in K$, 
the system $F(\o{x},H_d(\o{x}),\o{z})=0$ has 
an isolated solution which satisfies $||\o{x}|| \leq |a_d|^{d^c}-1/d$.
Let us choose $d_0$ so large that Proposition~\ref{induction} also
applies for $d \geq d_0$. 
The proof will be complete if we can show that the following claim is
true: for $d \geq d_0$ 
and any $\o{z} \in K$, the system 
$F(\o{x},H_d(\o{x}),\o{z})=0$ has an isolated solution in the closed ball
$\o{B}(0,R-1/d)$ where $R=|a_{d_0}|^{d_0^c}$.
The proof of this claim is a straightforward induction on $d$.
Indeed, the claim is true for $d_0$ by choice of $d_0$, and
one can go from $d$ to $d+1$  by Proposition~\ref{induction}.
\end{proof}
Although this is not really needed for the proof of our main result
we note that the same property does not hold only for rational balls,
but for all compact subsets of $\C^r$.
\begin{corollary} \label{compact}
Fix a polynomial map 
$F:\C^n \times \C^n \times \C^r \rightarrow \C^n$
with integer coefficients 
and a compact  $K \subseteq \C^r$ such that $F_{\o{z}}$ 
has a well-balanced zero for all $\o{z} \in K$.
There exists $R>0$ such that the following property 
holds if $d$ is sufficiently large: for any $\o{z} \in K$ the system 
$F(\o{x},H_d(\o{x}),\o{z})=0$ has an isolated 
 root in the closed ball $\o{B}(0,R)$.
\end{corollary}
\begin{proof}
By compactness of $K$ and Proposition~\ref{open}, $K$ can be covered
by a finite set $\{K_1,\ldots,K_p\}$ of rational balls such that
$F_{\o{\zeta}}$ has a well balanced zero for any $i \in
\{1,\ldots,p\}$ and any $\o{\zeta} \in K_i$. 
By Theorem~\ref{boundedroots} there exist $R_1,\ldots,R_p$ 
such that the following property 
holds for all sufficiently large $d$: for any $i \in \{1,\ldots,p\}$
and any $\o{z} \in K_i$, the system 
$F(\o{x},H_d(\o{x}),\o{z})=0$ has an isolated
 root in the closed ball $\o{B}(0,R_i)$.
Now set $R=\max(R_1,\ldots,R_p)$.
\end{proof}

\begin{proof}[of the main theorem]
We can write $G=F_{\o{z}}$  where $F: \C^n \times
\C^n \times \C^r \rightarrow \C^n$ is a polynomial map with integer
coefficients and $\o{z} \in \C^r$ is a tuple of parameters.
By Proposition~\ref{open} there exists a rational ball $K \subseteq \C^r$ 
containing $\o{z}$ such that $F_{\o{\zeta}}$ has a well balanced zero for
all $\o{\zeta}$ in $K$. 
By Theorem~\ref{boundedroots} there exists an increasing sequence
$(d_i)_{i \geq 0}$ of integers and a  sequence $(\o{x}_i)_{i \geq 0}$ 
of points of $\C^{n}$ such that $F(\o{x}_i,H_{d_i}(\o{x}_i),\o{z})=0$.
Moreover the sequence $(\o{x}_i)_{i \geq 0}$ 
remains  inside a fixed  compact ball $\o{B}(0,R)$.
Extracting a subsequence if necessary, we may therefore assume that
$\o{x}_i$ converges to a limit point $\o{a} \in \C^n$ as $i$ goes to
infinity. We conclude that $F(\o{a},H(\o{a}),\o{z})=0$ 
since $\lim_{d \rightarrow +\infty} F(\o{x},H_d(\o{x}),\o{z})
= F(\o{x},H(\o{x}),\o{z})$
uniformly with respect to $\o{x} \in \o{B}(0,R)$.
\end{proof}

\begin{remark}
We have only defined the notion of well balanced zero for a polynomial
map $G:\C^n \times \C^n \rightarrow \C^n$, but this notion clearly
makes sense if $G$ is an arbitrary analytic function.
The main theorem is no longer true in this more general context.
Indeed, set $G(x,y)=y-H(x)+1$. This function has well balanced zeros,
but there does not exist $a$ such that $G(a,H(a))=0$.
\end{remark}

\end{document}